\documentclass{article} 
\usepackage{amsmath}
\usepackage{amssymb}
\usepackage{amsthm}
\usepackage{amsfonts}
\def\ZZ#1#2{Z_{#2}^{#1}}
\def\EE#1#2{E_{#2}^{#1}}
\newtheorem{property}{Property}
\newtheorem{Lemma}{Lemma}
\newtheorem{definition}{Definition}
\newtheorem{theorem}{Theorem}

\newcommand\R{{\mathrm I}\! {\mathrm R}}

\def\mm{\mathfrak{m}}

\newcommand\Nor{\mathcal{N}}
\newcommand\Tran{\mathcal{T}}

\newcommand\Z{{\mathbb Z}}

\newcommand\N{{\mathbb N}}
\newcommand\Q{{\mathbb Q}}

\newcommand\p[1]{\frac{\partial}{\partial #1}}

\begin{document}
\bibliographystyle{alpha}
\title{Normal form theory and spectral sequences}
\author{Jan A. Sanders \\
Vrije Universiteit\\
Faculty of Sciences \\ 
Department of Mathematics\\ 
De Boelelaan 1081a\\
1081 HV Amsterdam\\
The Netherlands\\
E-mail: jansa@cs.vu.nl
}
\date{\today}
\abstract{
The concept of unique normal form is formulated in terms of a spectral sequence.
As an illustration of this technique some results of Baider and Churchill 
concerning
the normal form of the anharmonic oscillator are reproduced.
The aim of this paper is to show that spectral sequences give us a natural
framework in which to formulate normal form theory}
\section{Introduction}
The goal of this paper is to show that spectral sequences
give us a natural language in which to describe the process
of computing the unique normal form.
For a solid modern introduction to normal form theory and its historical development,
see \cite{MR1941477}.
Although this is in no way intended as a review of (unique) normal form theory,
let me give the reader some references to the literature which will at least
help to create the impression that this is a subject that still gets a lot of attention
nowadays: 
\cite{
MR1842050,
MR1865629,
MR2002b:37070,
MR2001a:37065,
MR2000j:34056,
MR2000b:34057,
MR99k:34080,
MR88d:58084,
MR87h:58167,
MR1834068,
MR1866913,
MR1847982,
MR2002b:34010,
MR2001i:37073,
MR2001g:37071,
MR1873273,
MR95i:58160,
MR90c:58150,
MR1875467,
MR1732242,
MR1855013,
MR2002e:34061,
MR2002d:37081,
MR2002d:34065,
MR2000h:37071,
MR99k:34079}.

Consider formal vectorfields at equilibrium and apply
formal transformations to them, assuming for the moment that the
linear part of the field is already in normal form (as in Jordan normal form).
A {\em unique normal form}, as it is called in the literature, is by no means unique.
But if two different procedures are used to compute the unique normal form,
with result \(N_1(v)\) and \(N_2(v)\), say, then one should have
\(N_1(v)=N_1(N_2(v))\) and vice versa.
If this is not the case, then surely at least one of the procedures
does not yield a unique normal form.
As observed by A. Baider, the uniqueness is more apparent
in the space of allowable transformations.
Once we know there are no transformations left to us,
the result is the unique normal form.

The process of computing normal forms consists of solving
the so-called homological equation. 
Much as this is supposed to remind one of homology,
this is never in any way used.
In this paper I will try and formulate normal form theory
in terms of cohomology, using the framework of spectral sequences.
The theory of spectral sequences is set up to do approximate calculations 
in filtered and graded differential modules, and therefore
is a likely candidate for a theoretical tool.
The technical problem is, that one would like to do all calculations
with the approximate normal form, since it is very important whether
certain coefficients in this normal form are invertible or not.
This has as a consequence that we do not have one differential operator,
but many. However, they converge in the filtration topology
and we can easily adapt the usual construction to suit our problem.

Having seen the theoretical part, the reader may still have absolutely no idea
what it means, so I have given as an illustration 
in section \ref{HO} one of the
few systems where the unique normal form can be completely understood,
namely the anharmonic oscillator, following the analysis
in \cite{MR90k:58146,MR90e:58135}.
For other examples, see \cite{MR92g:58113,MR93g:58131}.
In fact, the results as presented here are slightly more general than
Baider and Churchill's, since we allow coefficients in a local ring instead of a field.
To say that they are new would, however, be an overstatement.
The analysis does however provide us with a nice illustration
of the spectral sequence method. All proofs are reduced to elementary
calculations, and all the filtering arguments that complicate Baider's  proof
are already contained in the setup.
It is tempting to do same same for another class of equations,
namely those planar vectorfields with nilpotent linear part,
but one look at the length of the resulting analysis in \cite{MR93m:58101}
is convincing enough not to present this as an example, since it would
be another big paper, and no new results in it!
However, it is also clear that the basic construct in {\em loc.cit.} is
the Tic-Tac-Toe--construction, cf. section \ref{subs1}, and the whole spectral sequence approach
seems completely natural.
Since the whole analysis in this case relies on the judicious choice
of a second filtering, one probably needs to set up a context
with two filterings, and for each of these one constructs a coboundary
operator on a bicomplex \(\EE{q_1,q_2}{p_1,p_2}\).
The reader is encouraged to reformulate a familiar problem
in the language of spectral sequences, just to experience whether
this approach is as natural as claimed.

A similar approach has been used in \cite{MR58:31167,MR52:4332},\cite[Chapter 14]{MR86f:58018}
in the context of singularity theory.
Undoubtedly, from a higher point of view, one can view the present paper
as a simple corollary of the previous work by Arnol'd et.al..
Yet it took the author quite some time to figure this out.
So here we have an idea that has been around for a quarter of a century
and it has not been picked up by the normal form community,
despite the fact that in the meantime serious work was done to formulate
the theoretical basis in terms of filtered Lie algebras.
That means that now it is time to make some propaganda for the method and illustrate
its power.\footnote{Added in proof: The reader may also want to consult
\cite{Chu02} and \cite{Mur03}.}

{\bf Acknowledgment} The author would like to thank Jim Murdock for long e-mail
discussions and valuable questions and comments, and Andr{\'e} Vanderbauwhede
for listening to his explanations on their beach walks at Cala Gonone.
Thanks go to the Newton Institute
and the organizers of the Integrable Systems program
and the Netherlands Organization for Scientific Research (NWO)
for their financial support.
\section{Normal forms, first steps}
In the theory of normal forms of ordinary differential equations
at equilibrium, one studies
equations of the form 
\[
\dot{x}=v_0+v_1+\cdots,
\]
where \(v_i\) stands for the Taylor expansion of the vectorfield
of degree \(i+1\).
One computes a normal form by applying the formal transformation
\[
x=y+X_1+X_2+\cdots,
\]
and this leads to the so-called {\em homological equation}
\[
[ v_0 , X_1 ]= v_1-v_1^0,
\]
where \(v_1^0\) is the {\bf first order normal form} at the quadratic level.
Why this is called a homological equation is seldom explained
and this paper is written to provide an explanation of this terminology
and to define the so-called {\bf unique normal form}
in terms of spectral sequences.

Let us start by constructing a short complex as follows.
\(V\) is the space spanned by the given vector field
\(v=v_0+v_1+\cdots\), and let \(W\) be the space of vectorfields
starting with linear terms. Consider now the standard homology complex
(where \(W\) is seen as the representation space of the Lie algebra \(V\))
\[
0\leftarrow W \leftarrow V\otimes W \leftarrow \bigwedge^2 V\otimes W\cdots.
\]
Since \(\dim V=1\), \(\bigwedge^2 V=0\), and
one has the short sequence
\[
0\leftarrow W \leftarrow V\otimes W \leftarrow 0,
\]
where the only nonzero boundary map \(\partial\) is given by
\[
\partial(v\otimes w) = [v,w].
\]
Normal form theory of vectorfields is complicated by the fact
that the vectorfield with which we act is being changed at the same
time. Therefore we define a sequence of vectorfields
\(v^{j+1}=exp(ad(t^j))v^j\), with \(v^0=v\)
and one-dimensional spaces \(V^j=\langle v^j \rangle\), \(v^j\in W\).
The corresponding differentials can be composed as usual,
since one of the two will always be zero anyway.
\section{Spectral sequence}
\begin{quote}
"Spectral sequences were invented by Jean Leray, as a prisoner
of war during World War II, in order to compute the homology
(or cohomology) of a chain complex. They were made algebraic by Koszul
in 1945."
\cite{MR95f:18001}.
\end{quote}
We now formulate the general theory in terms of filtered Lie algebras,
and illustrate it by making comments on the interpretation in terms of vectorfields.
The whole construction will work equally well for Hamiltonian systems,
to mention one important class of examples for normal form theory.

Suppose that \(W\) is a filtered Lie algebra
\(W=W_0\supset W_1\supset W_2\supset\cdots\),
with \(\bigcap_{i=0}^\infty W_i=0\)
and \([W_i,W_j]\subset W_{i+j}\).
Consider the short complex
\[
0 \leftarrow W \stackrel{\partial_q}{\leftarrow} V^q \otimes W \leftarrow 0,
\]
with \(V^q=\langle v^q \rangle\), \(v^0=v\) and the higher 
\(v^q\) will be defined in the next section.
We assume here that \(v^{q+1}=v^q \pmod{W_{q+1}}\).
This corresponds to the fact that one does not perturb lower order terms
in the process of normalizing the vectorfield.
The filtration of \(W\) induces a filtration on \(\bigwedge^n V^q \otimes W\)
and we denote the filtration on the chain complex by
\[
F^q=F_0^q\supset\cdots \supset F_p^q\supset \cdots , \mbox{ with }
F_i^q=\sum_{n=0}^1 \bigwedge^n V^q \otimes W_i.
\]
For the normal form theory, we are not so much interested in the 
\(V^q\) part of these spaces, so instead we work with
\[
K=K_0\supset\cdots \supset K_p\supset \cdots , \mbox{ with }
K_i=W_i\oplus W_i=: \Nor_i^0\oplus \Tran_i^0.
\]
In many problems one can think of \(\Nor\) and \(\Tran\)
as essentially the same space, but if one considers for instance
problems with time-reversal symmetry, then one should take for \(\Tran\)
the elements that are invariant under the group action, and for \(\Nor\)
those that change sign (which is then compensated for by simultaneous time-reversal).

We write \(\partial_q(x,y)=([v^q,y],0)\), \(x\in \Nor_i^0, y\in \Tran_i^0\).
One should think of \(\Nor_i^0\) as the space where the normal form 
of the vectorfield lives, and \(\Tran_i^0\) is the space of transformations.

The reader may at this point wonder whether the fact that the boundary relations
are trivially satisfied will also trivialize the subsequent application of the
theory. The interesting thing is that it does not.
The interaction of the cohomology with the filtering is sufficiently
complicated to confuse anyone without a good organization.
It is the claim of this paper that spectral sequences provide
this organization.

The following is an almost  standard introduction, taken
from \cite{MR21:1583},
to spectral sequences,
which is reproduced here for the convenience of the reader.
For more information on spectral sequences, see \cite{MR2000m:55003,MR2002c:55027}.
The main difference with the usual approach is that we allow the \(\partial\)
to vary (by changing \(V^p\)). The assumption that \(v^j\) converges in the filtration topology
ensures that this does not damage the usual constructions.

We now consider an abstract homology complexes \[0\leftarrow K_p\stackrel{\partial_r}{\leftarrow}K_p\leftarrow 0,\]
with \(r=-1,0,1,\cdots\) and
\(\partial_{-1}=0\).
Now let \( \ZZ{r}{p} \) 
be the set of all \( x \in K_p \)
such that \( \partial_{r} x \in K_{p+r} \).
We assume the following property in some of our lemmas.
\begin{property}\label{dprop}
If \(x\in \ZZ{r}{p}\) or \(x\in \ZZ{r-1}{p}\), then \(\partial_rx-\partial_{r-1}x\in \ZZ{r-1}{p+r}\).
\end{property}
This property follows in the case of normal form theory from the fact that
the sequence of normal forms converges in the filtration topology,
that is to say, higher order calculations do not change lower order ones.
\begin{Lemma}
The complex defined above, with \(\partial_q(x,y)=([v^q,y],0)\), 
\(x\in \Nor_i^0, y\in \Tran_i^0\),
has Property \ref{dprop}.
\end{Lemma}
\begin{proof}
Since \(v^{r}=v^{r-1}+w_{r}\),
with \(w_{r}\in K_{r}\), \(\partial_{r}x=\partial_{r-1}x+([w_{r},x],0)\)
and \(([w_{r},x],0)\in \ZZ{r-1}{p+r}\).
\end{proof}
\begin{Lemma}
Let \(\{K_p,\partial_r\}_{p=0,r=-1}^\infty\) be homology complexes with Property \ref{dprop}.
Then \(  \ZZ{r-1}{p+1} \subset \ZZ{r}{p} \).
\end{Lemma}
\begin{proof}
If \( x \in \ZZ{r-1}{p+1} \) then \( x \in K_{p+1} \subset K_p \)
and \( \partial_{r-1} x \in K_{p+r} \). 
By Property \ref{dprop}, \(\partial_r x\in K_{p+r}\).
\end{proof}
\begin{Lemma}
Let \(\{K_p,\partial_r\}_{p=0,r=-1}^\infty\) be homology complexes with Property \ref{dprop}.
Then \( \partial_{r-1}  \ZZ{r-1}{p-r+1} \subset \ZZ{r}{p} \).
\end{Lemma}
\begin{proof}
If \( x \in \partial_{r-1}  \ZZ{r-1}{p-r+1} \)
then \( x = \partial_{r-1} y \) with \( y \in \ZZ{r-1}{p-r+1}\),
implying that \( y \in K_{p-r+1} \) and \( \partial_{r-1}  y \in K_p \).
It follows that \( x \in K_p \) and \( \partial_{r} x=0 \).
\end{proof}
\begin{definition}
Let \(\{K_p,\partial_r\}_{p=0,r=-1}^\infty\) be homology complexes with Property \ref{dprop}.
Then for \(r\geq 0\),
\begin{eqnarray}
\EE{r}{p} = \ZZ{r}{p}/(\partial_{r-1}  \ZZ{r-1}{p-r+1}+\ZZ{r-1}{p+1}) , 
\quad\EE{r}{\cdot} = \sum_p \EE{r}{p} .
\end{eqnarray}
\end{definition}
I would like to think as this as the natural definition of what the computation of
normal forms is all about. One divides out whatever can be transformed away,
the \(\partial_{r-1}  \ZZ{r-1}{p-r+1}\) term 
and the transformations that trivially give rise to higher order terms,
the \(\ZZ{r-1}{p+1}\) term.
The \(\ZZ{r}{p}\) itself sees to it that only those transformations are left that
do not perturb the previously computed lower order terms in normal form.
The next theorem then identifies these natural spaces in which the normal form 
and the allowable transformations live into cohomology spaces.
This identification greatly simplifies the calculation procedure
and allows us to use the familiar tools used in the analysis of bicomplexes.
\begin{theorem} \label{AppendixTheo3}
Let \(\{K_p,\partial_p\}_{p=0,r=-1}^\infty\) be homology complexes with Property \ref{dprop}.
Then there exists on the graded module \(\EE{r}{\cdot}\) a differential \(d_\cdot^r  \)
such that \( H(\EE{r}{\cdot}) \) is canonically 
isomorphic to \(\EE{r+1}{\cdot}\), \(r\geq 0\).
\end{theorem}
\begin{proof}
Here only the definition of \(d_\cdot^r  \) is given.
For the full proof, see Appendix \ref{App}.
The differential \( \partial_{r}  \) maps \( \ZZ{r}{p} \) into \( \ZZ{r}{p+r} \)
and \( \partial_{r-1}  \ZZ{r-1}{p-r+1}+\ZZ{r-1}{p+1}\) into \( \partial_{r} \ZZ{r-1}{p+1}\).
Let \(x\in \partial_{r} \ZZ{r-1}{p+1}\). Then there is a \(y\in \ZZ{r-1}{p+1}\) such that
\(x=\partial_r y \in \partial_{r-1}\ZZ{r-1}{p+1}+\ZZ{r-1}{p+r+1}\), using Property \ref{dprop}.
Since
\begin{eqnarray}
\EE{r}{p+r} = \ZZ{r}{p+r}/(\partial_{r-1}  \ZZ{r-1}{p+1}+\ZZ{r-1}{p+r+1}) ,
\end{eqnarray}
one sees that \( \partial_{r}  \)  induces a 
map \( d_p^r  : \EE{r}{p} \rightarrow \EE{r}{p+r} \).
\end{proof}
\section{Normal form theory}
The direct sum decomposition \(K_{p}=\Nor_p^{0}\oplus \Tran_p^{0} \)
induces a analogous decomposition \(\EE{r}{p}=\Nor_p^{r}\oplus \Tran_p^{r} \).
We will see from the normal form calculations that this decomposition
codifies the form of the normal terms at level \(r\) and degree \(p+1\)
in \(\Nor_p^{r}\) and the terms we can still use in the transformation in 
\(\Tran_p^{r}\).
Let us see whether this whole construction makes some sense
in terms of classical normal form theory.
Take \(v^0=v\), the vectorfield we start with.
We construct \(\EE{1}{p}\). By definition
\begin{eqnarray}
\EE{1}{p} = \ZZ{1}{p}/(\partial_0  \ZZ{0}{p}+\ZZ{0}{p+1}) .
\end{eqnarray}
We know that \(\ZZ{1}{p}=\{x\in K_p | \partial_0 x \in K_{p+1}\}\).
We write for each \(x\in\ZZ{1}{p}\), \(x=(n, t)\), with \(t\) such that \([v^0,t]\in W_{p+1}\).
Dividing by \(\partial_0  \ZZ{0}{p}\) means we put \(n\) in first order normal form,
that is, \(n\in\Nor_p^{1}\),
since \(\partial_0  \ZZ{0}{p}\) is exactly the image of \(ad(v^0)\) restricted to the 
\(p+1\)-degree
terms both in the source and the object space.
Thus \(\EE{1}{p}\) can be identified with pairs \((n,t)\),
with \(n\) in first order normal form (with respect to the linear part of
the vector field \(v^0\)) and \(t\in \Tran_p^1\) such that it does not change the \(p+1\)-degree terms
when we use it to transform the equation.
This means that the lowest order term in \(t\) commutes
with the linear term of the vectorfield.

In the following we use the subindex to indicate the graded part
induced by the filtration, as in \(v_p\in W_p/W_{p+1}=:G_p\).
In the context of formal power series vectorfields these are
just the homogeneous parts of degree \(p+1\).

We now write \(v_1^0=v_1^1+[v_0^0,t_1^0]\), with
\([(v_1^1,0)]\in \EE{1}{1}\) (where the \([\cdot]\) denote equivalence classes,
not Lie brackets!) and \([(0,t_1^0)]\in \EE{0}{1}\).
The decomposition of \(G_1\) implicit in the equality amounts to the identification
\(G_1\approx (G_1/im(ad(v_0^0)|G_1))\oplus \ker(ad(v_0^0)|G_1)\),
which in turn is equivalent to the choice of a complement for \(im(ad(v_0^0)\) in \(G_1\).
Choosing such a complement is standard in normal form theory
(which is not to say there is a standard choice!) and we now see a homological formulation.
We define
\[
v^1=exp(ad(t_1^0))v^0.
\]
We see that
\begin{eqnarray*}
v^1&=&v_0^0+v_1^0+ad(t_1^0)v_0^0 \pmod{W_2}
\\& =&v_0^0+v_1^1+[v_0^0,t_1^0]+[t_1^0,v_0^0]\pmod{W_2} 
\\&=&v_0^0+ v_1^1\pmod{W_2}.
\end{eqnarray*}
We are now ready to go to round two.
We write \(v_2^1=v_2^2+[v_0^0,t_2^0]+[v_0^0+v_1^1,t_1^1]\),
with \(v_2^2\in \Nor_2^2\), \(t_2^0\in \Tran_2^0=W_2\) and \(t_1^1\in \Tran_1^1\),
with \(\EE{2}{2}=\Nor_2^2\oplus \Tran_2^2\).
Let us now compute \(d_p^1 (v_p^1,t_p^1)=(\pi_{p+1}^1[v_1^0,t_p^1],0)\),
where \(\pi_{p+1}^1 \) is the projection on \(\Nor_{p+1}^1\) and \(d_p^1\)
is defined in the proof of Theorem \ref{AppendixTheo3}.
So this gives us the contribution of the transformation
to the normal form of the terms of one degree higher,
and we see that the general spectral sequence approach
gives us exactly those terms that we have always been computing.
This means that this is the natural language in which to formulate 
the theory.

Then we define
\[
v^2=exp(ad(t_1^1+t_2^0))v^1.
\]
In general we write
\begin{eqnarray*}
(v_{j+1}^j,0)&=&(v_{j+1}^{j+1},0)+\sum_{k=0}^{j} ([\sum_{i=0}^{k} v_i^i , t_{j+1-k}^k],0)
\\&=&(v_{j+1}^{j+1},0)+\sum_{k=0}^{j} \partial_{k} (0, t_{j+1-k}^k)
\\&=&
(v_{j+1}^{j+1},0)+\partial_{j}(0,\sum_{k=0}^{j} t_{j+1-k}^k) \pmod{K_{j+2}}
\end{eqnarray*}
and transform with
\[
v^{j+1}=exp(ad(\sum_{k=0}^{j} t_{j+1-k}^k )) v^j.
\]
Here one has \([(v_{j+1}^{j+1},0)]\in Z_{j+1}^{j+1}\),
and \([\partial_{j}(0,\sum_{k=0}^{j} t_{j+1-k}^k)]\in \partial_{j} Z_1^j +Z_{j+2}^{j}\).
Thus one can view \((v_{j+1}^{j+1},0)\) as an element in \(\EE{j+1}{j+1}\).

One sees that the definitions of the spaces allow one to find all
these elements in principle, but one has to be aware of possible problems
if the transformation space is not locally finite,
that is the factor spaces are not finite dimensional.
This is for instance the case if one considers
equations of the type
\[
\dot{x}=\epsilon f(t,x).
\]
If \(f\) is periodic in \(t\),
one can use averaging to solve the homological equation \cite{MR95g:58209}\nocite{MR94m:00024}.
The filtering is given here by the powers of \(\epsilon\).
If they are finite dimensional, everything is just linear algebra.
Most of normal form theory is about doing the linear algebra in an effective
way 
by using the 
spectral information of the linear part of the vectorfield,
but this only works well for the first order normal form,
since the Lie algebra involved in the higher order normal form
is not reductive (Try to imbed an element in \(W_1\) in an \(sl(2,\R)\)).
This is treated in \cite{MR95g:58209}.

Continuing in this fashion, we decrease the normal form space and the space of transformations
we can use, until in the end we have the unique normal form space, and the transformations
that are left commute with the normal form of \(v\), that is, they are 
conjugate to symmetries of \(v\).
\begin{theorem}
Suppose \(\EE{\infty}{\cdot}\) exists and let \(v^\infty\) be the final normal form
of \(v\). 
Then we can identify
\( v_p^\infty\) with an element in \(\EE{p}{p}=\EE{\infty}{p}\)
of the form \((v_p^\infty,0)\).
Furthermore we can identify any symmetry \(s\) of \(v\) with
a symmetry \(s^\infty\) of \(v^\infty\) and \((0,s_p^\infty)\in \EE{\infty}{p}\).
\end{theorem}
\section{The anharmonic oscillator}\label{HO}
\begin{quote}
"[] the most powerful method of computing homology groups uses spectral sequences.
When I was a graduate student, I always wanted to say, nonchalantly,
that such and such is true "by the usual spectral sequence argument,"
but I never had the nerve"
\cite{Ro96}.
\end{quote}
Let us, just to get used to the notation, treat the simplest
normal form problem we can think of, the anharmonic oscillator.
The results we obtain were obtained first in \cite{MR90k:58146,MR90e:58135}.

We will take our coefficients from a local ring \(R\) containing \(\Q\). 
\footnote{One can think for instance of formal power series in a deformation parameter \(\lambda\),
which is the typical situation in bifurcation problems.
Then a term \(\lambda x^2\p{x}\) has coefficient \(\lambda\) which
is neither zero nor invertible, since \(\frac{1}{\lambda}\) is not a formal
power series.}
Then the noninvertible
elements are in the maximal ideal, say \(\mm \),
and although subsequent computations are going to affect
terms that we already consider as fixed in the normal form calculation, 
they will not affect their equivalence class in the residue field \(R/\mm \).
So the convergence of the spectral sequence is with respect to
the residue field. The actual normal form will contain formal power
series which converge in the \(\mm \)-adic topology.
In \cite{MR90e:58135} it is assumed there is no maximal ideal,
and \(R=\R\).
When in the sequel we say that something is in the kernel of
a coboundary operator, this means that the result has its coefficients in \(\mm \).
When we compute the image then this is done first in the residue field
to check invertibility, and then extended to the whole of \(R\).
This gives us more accurate information than simply listing the normal form
with coefficients in a field,
since it allows for terms which have nonzero coefficients, through which
we do not want to divide, either because they are very small
or because they contain a deformation parameter in such a way that the coefficient
is zero for one or more values of this parameter.
\begin{definition}
We define the Hilbert-Poincar{\'e} series of \(\EE{r}{\cdot}\) as
\[
P[\EE{r}{\cdot}](t)=\sum_{p=0}^\infty (\dim \Nor_p^r-\dim \Tran_p^r)t^p.
\]
\end{definition}
If it exists, we call \(I(\EE{r}{\cdot})=P[\EE{r}{\cdot}](1)\) the {\bf index}
of the spectral sequence at \(r\).
In the anharmonic oscillator problem, \(P[\EE{0}{\cdot}](t)=0\).
Let,  with \(k\geq -1, l\geq 0, q\in \Z/4\),
\(A_{k+l}^{k-l,q}=i^q(x^{k+1}y^l\p{x} +i^{2q} x^l y^{k+1} \p{y})\).
Since \(A_k^{l,q+2}=-A_k^{l,q}\),  a basis is given by 
\(\langle A_k^{l,q} \rangle_{k=-1,\cdots,l=0,\cdots,q=0,1}\),
but we have to compute in \(\Z/4\).
The commutation relation is
\begin{eqnarray*}
\lefteqn{[A_{k+l}^{k-l,p},A_{m+n}^{m-n,q}]=}&&\\
&=&(m-k)A_{k+m+l+n}^{k-l+m-n,p+q}
+nA_{k+m+l+n}^{m-n-(k-l),q-p}-lA_{k+m+l+n}^{k-l-(m-n),p-q}.
\end{eqnarray*}
Then the anharmonic oscillator is of the form
\[
v=A_{0}^{0,1}+\sum_{q=0}^1\sum_{k+l=1}^\infty \alpha_{k+l}^{k-l,q} A_{k+l}^{k-l,q}, \alpha_k^l \in R.
\]
Since
\[
[A_{0}^{0,1},A_{k+l}^{k-l,q}]=
(k-l)A_{k+l}^{k-l,q+1}
\]
we see that the kernel of \(ad(A_{0}^{0,1})\) consists of those
\(A_{k+l}^{k-l,q}\) with \(k=l\)
and the image of those with \(k\neq l\).
We are now in a position to compute \(\EE{1}{p}\).
We have by definition that \(\ZZ{0}{p}=K_p=W_p\oplus W_p\).
Then
\begin{eqnarray*}
\ZZ{1}{p}&=&\{x\in K_p|\partial_1x\in K_{p+1}\}
\\&=&
\{(x,y)\in K_p|[v_0^0,y_p]=0 \}
\\&=& W_p\oplus \ker \ ad(v_0^0)|_{W_p}.
\end{eqnarray*}
Thus, 
\[
\EE{1}{p}=\ZZ{1}{p}/(\partial_0\ZZ{0}{p}+\ZZ{0}{p+1})=\ker \ ad(v_0^0)|_{G_p}\oplus \ker \ ad(v_0^0)|_{G_p},
\]
since \(W_p=im\ ad(v_0^0)|_{W{p}}\oplus \ker\ ad(v_0^0) |_{W{p}}\) and
\(G_p=im\ ad(v_0^0)|_{G{p}}\oplus \ker\ ad(v_0^0) |_{G{p}}\),
due to the semisimplicity of \(A_0^{0,1}\).
It follows that \(P[\EE{1}{\cdot}](t)=0\).
In general we have
\[
\Tran_{2p}^1=\Nor_{2p}^1={\langle A_{2p}^{0,0}, A_{2p}^{0,1} \rangle}_{R},
\]
and \(\Tran_{2p+1}^1=\Nor_{2p+1}^1=0\). 
One has the following commutation relations
\begin{eqnarray*}
[A_{2k}^{0,p},A_{2m}^{0,q}]
=(m-k)A_{2k+2m}^{0,p+q}
+mA_{2k+2m}^{0,q-p}
-kA_{2k+2m}^{0,p-q}.
\end{eqnarray*}
For later use we write out the three different cases:
\begin{eqnarray*}
&& [A_{2k}^{0,0},A_{2m}^{0,0}]
=2(m-k)A_{2k+2m}^{0,0},
\\&&
[A_{2k}^{0,0},A_{2m}^{0,1}]
 = 2mA_{2k+2m}^{0,1},
\\&&
[A_{2k}^{0,1},A_{2m}^{0,1}]
=0.
\end{eqnarray*}
It follows that \(\mathcal{A}=\langle A_{2m}^{0,1}\rangle_{m\in\N} \oplus \langle A_{2m}^{0,1}\rangle_{m\in\N} \)
is an abelian Lie algebra ideal in \(\EE{1}{\cdot}\),
which itself is a \(\N\times \Z/2\)-graded Lie algebra.
We can consider \(\EE{1}{\cdot}\) as a central extension
of \(\EE{1}{\cdot,0}\) with \(\EE{1}{\cdot,1}\).

We now continue our normal form calculations until we hit a term
\(v_{2r}^{2r}=\beta_{2r}^0 A_{2r}^{0,0}+\beta_{2r}^1 A_{2r}^{0,1}\)
with either \(\beta_{2r}^0\) or \(\beta_{2r}^1\) invertible.
We have \(\EE{2r}{\cdot}=\EE{1}{\cdot}\).
We see that \(d_{2p}^{2r}\) is now nonzero, at least it is not
zero by the previous argument, since it maps the even spaces
on themselves.
A general element in \(\Tran_{2p}^{2r} \) is given by \begin{eqnarray*}
t_{2p}^{2r}&=&\sum_{q=0}^1 \gamma_{2p}^q A_{2p}^{0,q}.
\end{eqnarray*}
We have, with \(p>r\),
\begin{eqnarray*}
\lefteqn{d_{2p}^{2r}(0,t_{2p}^{2r})=}&&
\\&=&
(\beta_{2r}^0 \gamma_{2p}^0[A_{2r}^{0,0},A_{2p}^{0,0}] 
+\beta_{2r}^0 \gamma_{2p}^1[A_{2r}^{0,0},A_{2p}^{0,1}] 
+\beta_{2r}^1 \gamma_{2p}^0[A_{2r}^{0,1},A_{2p}^{0,0}] ,0)
\\&=&
(2(p-r) \beta_{2r}^0 \gamma_{2p}^0A_{2p+2r}^{0,0}
+2p \beta_{2r}^0 \gamma_{2p}^1A_{2p+2r}^{0,1}
-2r\beta_{2r}^1 \gamma_{2p}^0A_{2p+2r}^{0,1},0).
\end{eqnarray*}
We view this as a map from the coefficients at \(G_{2p}\) to those at \(G_{2r+2p}\)
with matrix representation
\[
\left(
\begin{array}{cc}
2(p-r)\beta_{2r}^0 &0 \\ -2r \beta_{2r}^1  & 2p \beta_{2r}^0
\end{array}\right)
\left(
\begin{array}{c}
\gamma_{2p}^0 \\ \gamma_{2p}^1
\end{array}\right)
\]
and we see that for \(0< p\neq r \) the map is surjective if \(\beta_{2r}^0\) is invertible;
if it is not, it has a one-dimensional image since we assume
that in this case \( \beta_{2r}^1\) is invertible.
\subsection{Case \(\mathcal{A}^r\): \(\beta_{2r}^0\) is invertible.}
In this subsection we assume that \(\beta_{2r}^0\) is invertible.
The following analysis is equivalent to the one in \cite[Theorem 4.11]{MR90e:58135},
case (3), \(j=r\), if \(\beta_{2r}^1=0\).
For \(\beta_{2r}^1\neq 0\), see section \ref{subs}.

We have already shown that \(im\ d_{2p}^{2r}=\langle A_{2p+2r}^{0,0},A_{2p+2r}^{0,1}\rangle
\oplus 0\) for \(0<p\neq r\) and \(im\ d_{2r}^{2r}=\langle A_{4r}^{0,1}\rangle \oplus 0\).
Furthermore \(\ker \ d_{2r}^{2r}=G_2 \oplus \langle v_{2r}^{2r}\rangle\)
and \(\ker \ d_{2p}^{2r}=G_2 \oplus 0\) for \(0<p\neq r\).
We are now in a position to compute \(\EE{2r+1}{2p}=H^{2p}(\EE{2r}{\cdot})\).
First of all, \(\EE{2r+1}{4r}=H^{4r}(\EE{2r}{\cdot})=\ker \ d_{4r}^{2r} / im \ d_{2r}^{2r}
=( \langle A_{4r}^{0,0}, A_{4r}^{0,1}\rangle \oplus 0 )/(\langle A_{4r}^{0,1}\rangle \oplus 0)=\langle A_{4r}^{0,0}\rangle\oplus 0
\).
Then for \(p>r\) we find
\(\EE{2r+1}{2p} = H^{2p}(\EE{2r}{\cdot})=\ker \ d_{2p}^{2r} / im \ d_{2p-2r}^{2r}
=(\langle A_{2p}^{0,0}, A_{2p}^{0,1}\rangle \oplus 0)/(\langle A_{2p}^{0,0}, A_{2p}^{0,1}\rangle \oplus 0)=0\), while for \(0<p<r \) we find
\(\EE{2r+1}{2p} = H^{2p}(\EE{2r}{\cdot})=\ker \ d_{2p}^{2r} 
=(\langle A_{2p}^{0,0}, A_{2p}^{0,1}\rangle \oplus 0)\).
Obviously, \(\EE{2r+1}{2p+1}=0\).
One has
\[
\EE{\infty}{2p}=\EE{2r+1}{2p}=\left\{ \begin{array}{cccccl}
A_{2p}^{0,0}&A_{2p}^{0,1}&\oplus&   
A_{2p}^{0,0}&A_{2p}^{0,1}&p
\\
\\
 \mm &R\setminus\mm  &&R&R&  0 
\\
\mm &\mm & & 0&0& 1,\cdots,r-1 
\\
R\setminus\mm &R& & R&0& r 
\\
0&0&&0&0&  r+1,\cdots,2r-1
\\
 R &0& & 0&0& 2r
\\
0&0&& 0 &0&  2r+1,\cdots .
\end{array}\right. 
\]
We see that
\[
P^r[\EE{\infty}{\cdot}](t)=\sum_{i=1}^{r-1}2t^{2i}+t^{2r}+t^{4r}\]
and \(I(\EE{\infty}{\cdot})=2r\).
The codimension of the sequence, which we obtain by looking at the dimension
of the space with coefficients in \(\mm \), is \(2r-1\).
We can reconstruct the normal form out of this result.
Here \(c_{2p}\in \mm \) at position \(A_{2p}^{0,\cdot}\) means that the coefficient
of \(A_{2p}^{0,\cdot}\) in \(c_{2p}\) cannot be invertible.
And \(c_{2p}\in R\setminus\mm \)  means that it should be invertible.
While \(c_{2p}\in R\) indicates that the coefficient could be anything
in \(R\).
By ignoring the \(\mm A_{2p}^{0,\cdot}\) terms we obtain the results
in \cite{MR90e:58135}.
The \(R\setminus\mm \)-terms indicate the organizing center of the corresponding
bifurcation problem.

Since we have no more effective transformations at our disposal, all cohomology
after this will be trivial and we have reached the end of our spectral sequence
calculation.

\subsection{Case \(\mathcal{A}_r\): \(\beta_{2r}^0\) is not invertible,
but \( \beta_{2r}^1\) is}\label{subs}
The following analysis is equivalent to the one in \cite[Theorem 4.11]{MR90e:58135},
case (4), \(k=r, l=q\).
Since
\begin{eqnarray*}
\lefteqn{d_{2p}^{2r}(0,t_{2p})=}&&
\\&=&
2(p-r) \beta_{2r}^0 \gamma_{2p}^0A_{2p+2r}^{0,0}
+2p \beta_{2r}^0 \gamma_{2p}^1A_{2p+2r}^{0,1}
-2r\beta_{2r}^1 \gamma_{2p}^0A_{2p+2r}^{0,1}
\end{eqnarray*}
we can remove all terms \(A_{2p+2r}^{0,1}\) for \(p>0\)
by taking  \(\gamma_{2p}^1=0\).
This only contributes terms in \(\mm  A_{2p+2r}^{0,0}\).
We obtain
\[
\EE{2r+1}{2p}=\left\{ \begin{array}{cccccl}
A_{2p}^{0,0}&A_{2p}^{0,1}& \oplus&
A_{2p}^{0,0}&A_{2p}^{0,1}&p
\\
\\
\mm &R\setminus\mm & &R&R& 0
\\
\mm &\mm && 0&R & 1,\cdots,r-1
\\
\mm &R\setminus\mm && 0&R &  r
\\
 R&0 && 0&R & r+1,\cdots .
\end{array}\right.
\]
We see that
\[
P_r[\EE{2r+1}{\cdot}](t)=\sum_{i=1}^{r}t^{2i}\]
and \(I(\EE{2r+1}{\cdot})=r\).
The codimension is \(2r\).
\subsubsection{Case \(\mathcal{A}_r^q\): \( \beta_{2q}^0\) is invertible}\label{subs1}
We now continue our normal form calculation until at some point we hit on a term
\[
\beta_{2q}^0 A_{2q}^{0,0}
\]
with \(\beta_{2q}^0 \) invertible.
The following argument is basically 
the Tic-Tac-Toe
Lemma \cite[Proposition 12.1]{MR83i:57016}
and this was a strong motivation to consider
spectral sequences as  a framework for normal form theory.
The idea is to add the \(\Z/2\)-grading to our considerations.
We view \(ad(A_0^{0,1}+\beta_{2r}^1 A_{2r}^{0,1})\) as one coboundary operator
\(d_{\cdot}^{2r,1}\) and \(ad(\beta_{2q}^{0} A_{2q}^{0,0})\) as another,
\(d_{\cdot}^{2q,0}\). Both operators act completely homogeneous
with respect to the gradings induced by the filtering 
and allow us to consider the bicomplex spanned by
\(\EE{1}{\cdot,0}\) and \(\EE{1}{\cdot,1}\),
 where \(\EE{1}{2p,0}=\langle A_{2p}^{0,0} \rangle \oplus \langle A_{2p}^{0,0} \rangle\) and
 \(\EE{1}{2p,1}=\langle A_{2p}^{0,1}\rangle\oplus \langle A_{2p}^{0,1}\rangle\).
Since \([ A_{2q}^{0,0} , A_{2p}^{0,0}]=2(p-q)A_{2p+2q}^{0,0}\) and 
\([ A_{2q}^{0,0} , A_{2p}^{0,1}]=2pA_{2p+2q}^{0,1}\)
we see that the only nontrivial \(d^{2q,0}\)-cohomology is
\(H_{d^{2q,0}}^{4q} (\EE{1}{\cdot,0})\oplus H_{d^{2q,0}}^{2q} (\EE{1}{\cdot,0})\).

To compute the image of \(d_{\cdot}^{2r,1}+d_{\cdot}^{2q,0}\)
we start with the \(\EE{2r+1}{2s}\)-term. Take \(t_s^{(1)}=A_{2s}^{0,1}\).
Then \(d_{2s}^{2q,0}A_{2s}^{0,1}\in \EE{2r+1}{2q+2s,1}\),
that is, we can write \(d_{2s}^{2q,0}A_{2s}^{0,1}+d_{2q+2s-2r}^{2r,1}t_s^{(2)}=0\), with \(t_s^{(2)}\in \EE{1}{2q+2s-2r}\).
If we now compute \((d_{\cdot}^{2r,1}+d_{\cdot}^{2q,0})(t_s^{(1)}
+t_s^{(2)})\), we obtain
\[
(d_{\cdot}^{2r,1}+d_{\cdot}^{2q,0})(t_s^{(1)}
+t_s^{(2)})=d_{2q+2s-2r}^{2q,0}t_s^{(2)}.
\]
Looking at the \(d^{2q,0}\)-cohomology we see that this gives us a nonzero result
under the condition \(0<s\neq r\).

The image of \(d_{\cdot}^{2q}=d_{\cdot}^{2r,1}+d_{\cdot}^{2q,0}\) 
in \(\EE{1}{\cdot,0}\) is spanned by
\[ \prod_{\stackrel{j=1}{j\neq r}}^{\infty}\langle A_{4q-2r+2j}^{0,0}\rangle \oplus 0\]
and in \(\EE{1}{\cdot,1}\) by 
\( \prod_{j=1}^{\infty}\langle A_{2r+2j}^{0,1}\rangle \oplus 0\).
The kernel is spanned by \(A_{2r}^{0,1}\), that is, an element with this as its
lowest order term. This, of course, is the equation itself. 
Thus
\[
\EE{\infty}{2p}=\EE{2q+1}{2p}=\left\{ \begin{array}{cccccl}
A_{2p}^{0,0}&A_{2p}^{0,1}& \oplus&
A_{2p}^{0,0}&A_{2p}^{0,1}&p
\\
\\
\mm &R\setminus\mm & &R&R& 0
\\
\mm &\mm && 0&0 & 1,\cdots,r-1
\\
\mm &R\setminus\mm && 0&R &  r 
\\
 \mm &0 && 0&0 & r+1,\cdots,q -1
\\
R\setminus\mm  &0&& 0&0 & q 
\\
 R&0 && 0&0 & q+1,\cdots,2q-r 
\\
0&0&& 0 &0& 2q-r+1,\cdots,2q-1
\\
R &0&& 0&0 & 2q 
\\
0&0&& 0&0 & 2q+1,\cdots
\end{array}\right.
\]
We see that
\[
P_r^q[\EE{\infty}{\cdot}](t)=\sum_{i=1}^{r-1}2t^{2i}+t^{2r}+\sum_{i=r+1}^{q-1} t^{2i}
+\sum_{i=q}^{2q-r} t^{2i}+t^{4q}\]
and \(I(\EE{\infty}{\cdot})=2q\).
The codimension is \(r+q-1\).
This is the final result, since there is nothing useful to do for the
\(A_{2r}^{0,1}\) term in \(\EE{2q+1}{2r}\).
The \(A_{0}^{0,0}\)-term may be used to scale one of the coefficients
in \(R\setminus\mm \) to unity.
\subsubsection{Case \(\mathcal{A}_r^\infty\): no \( \beta_{2q}^0\) is invertible}\label{subs2}
The following analysis is equivalent to the one in \cite[Theorem 4.11]{MR90e:58135},
case (2), \(k=r\).

Since we can eliminate all terms of type \(A_{2p}^{0,1}\),
and we find no terms of type \(A_{2p}^{0,0}\) with invertible coefficients,
we can draw the conclusion that the cohomology is spanned by the
\(A_{2p}^{0,0}\), but does not show up in the normal form.
\[
\EE{\infty}{2p}=\EE{2r+1}{2p}=\left\{ \begin{array}{cccccl}
A_{2p}^{0,0}&A_{2p}^{0,1} &\oplus& 
A_{2p}^{0,0}&A_{2p}^{0,1} &p 
\\
\\
\mm &R\setminus\mm  &&R&R& 0 
\\
\mm &\mm  && 0&R& 1,\cdots,r-1
\\
\mm &R\setminus\mm  && 0&R &  r 
\\
\mm  &0&& 0&R&  r+1,\cdots
\end{array}\right.
\]
We see that
\[
P_r^\infty[\EE{\infty}{\cdot}](t)=\sum_{i=1}^{r}t^{2i}
\]
and \(I(\EE{\infty}{\cdot})=r\).
The codimension is infinite.
Scaling the coefficient of \(A_{2r}^{0,1}\) to unity
uses up the action of \(A_{0}^{0,0}\).
Although we still have some freedom in our choice of transformation,
this freedom cannot effectively be used, so it remains in the final result.
We summarize the index results as follows.

The index of \(\mathcal{A}_r^q\)
is \(2q\) if \(q\in\N\) and \(r\) otherwise.
\subsection{The \(\mm \)-adic approach}
So far we have done all computations
modulo \(\mm\). One can now continue doing the same thing, but
now on the \(\mm\) level, and so on. The result will be a finite sequence
of \(\mm^p \mathcal{A}_{r_p}^{q_p}\) describing exactly what remains.
Here the lower index can be either empty, a natural number or infinity
and the upper index can be a (bigger) natural number or infinity.
The generating function will be
\[
P[\EE{\infty}{\cdot}](t)=\sum_p u^p P_{r_p}^{q_p}[\EE{\infty}{\cdot}](t),
\]
with \(u^p\) standing for an element in \(\mm^p\setminus \mm^{p+1}\).
\appendix
\section{The spectral sequence theorem}\label{App}
\begin{theorem} 
There exists on the graded module \(\EE{r}{\cdot}\) a differential \(d_\cdot^r  \)
such that \( H(\EE{r}{\cdot}) \) is canonically 
isomorphic to \(\EE{r+1}{\cdot}\), \(r\geq 0\).
\end{theorem}
\begin{proof}
We follow \cite{MR21:1583} with modifications to allow for the 
converging boundary operators.
The differential \( \partial_{r}  \) maps \( \ZZ{r}{p} \) into \( \ZZ{r}{p+r} \)
and \( \partial_{r-1}  \ZZ{r-1}{p-r+1}+\ZZ{r-1}{p+1}\) into \( \partial_{r} \ZZ{r-1}{p+1}\).
Let \(x\in \partial_{r} \ZZ{r-1}{p+1}\). Then there is a \(y\in \ZZ{r-1}{p+1}\) such that
\(x=\partial_r y \in \partial_{r-1}\ZZ{r-1}{p+1}+\ZZ{r-1}{p+r+1}\), using Property \ref{dprop}.
Since
\begin{eqnarray}
\EE{r}{p+r} = \ZZ{r}{p+r}/(\partial_{r-1}  \ZZ{r-1}{p+1}+\ZZ{r-1}{p+r+1}) 
\end{eqnarray}
one sees that \( \partial_{r}  \)  induces a 
map \( d_p^r  : \EE{r}{p} \rightarrow \EE{r}{p+r} \).
For \(x\in \ZZ{r}{p} \) to define a cocycle of degree \(p\) on \(\EE{r}{\cdot}\)
it is necessary and sufficient that \( \partial_{r} x \in \partial_{r-1}  \ZZ{r-1}{p+1}+\ZZ{r-1}{p+r+1}\),
i.e. (again using Property \ref{dprop}) \( \partial_{r} x = \partial_{r} y + z \) with \( y \in \ZZ{r-1}{p+1} \) and \( z \in \ZZ{r-1}{p+r+1} \).
Putting \( u = x - y \in \ZZ{r}{p}+\ZZ{r-1}{p+1} \subset K_p \), with
\( \partial_{r} u=\partial_{r} x-\partial_{r} y = z \in K_{p+r+1} \), one has
\( u \in \ZZ{r+1}{p}\),
since \(\partial_{r+1}u -\partial_r u \in K_{p+r+1}\).
In other words, \( x \in \ZZ{r-1}{p+1}+\ZZ{r+1}{p} \).
It follows that the \(p\)-cocycles are given by
\begin{eqnarray}
Z^{p}(\EE{r}{\cdot})=(\ZZ{r+1}{p}+\ZZ{r-1}{p+1})/(\partial_{r-1} \ZZ{r-1}{p-r+1}+\ZZ{r-1}{p+1}) .
\end{eqnarray}
The space of \(p\)-coboundaries \( B^p(\EE{r}{\cdot}) \) 
consists of elements of \( \partial_{r}  \ZZ{r}{p-r} \)
and one has
\begin{eqnarray}
B^p(\EE{r}{\cdot})=(\partial_{r} \ZZ{r}{p-r}+\ZZ{r-1}{p+1})/(\partial_{r-1} \ZZ{r-1}{p-r+1}+\ZZ{r-1}{p+1}) .
\end{eqnarray}
It follows that
\begin{eqnarray}
H^p(\EE{r}{\cdot})&=&
(\ZZ{r+1}{p}+\ZZ{r-1}{p+1})/(\partial_{r} \ZZ{r}{p-r}+\ZZ{r-1}{p+1})\nonumber\\
&=& \ZZ{r+1}{p}/(\ZZ{r+1}{p}\cap(\partial_{r} \ZZ{r}{p-r}+\ZZ{r-1}{p+1})).
\end{eqnarray}
We now first prove that 
\( \ZZ{r+1}{p}\cap \ZZ{r-1}{p+1} = \ZZ{r}{p+1}\).
Let \( x \in \ZZ{r+1}{p} \cap \ZZ{r-1}{p+1}\).
Then \( x \in K_{p+1} \) and \( \partial_{r+1} x \in K_{p+r+1} \).
Thus \(\partial_r x\in K_{p+r+1}\) according to Property \ref{dprop}.
This implies \( x \in \ZZ{r}{p+1} \).
On the other hand, if \( x \in  \ZZ{r}{p+1}\)
we have \( x \in K_{p+1} \subset K_p \)
and \( \partial_{r} x \in K_{p+r+1} \subset K_{p+r} \).
Thus \( x \in K_p \) and \( \partial_{r} x \in K_{p+r+1}\).
Again it follows that \(\partial_{r+1} x \in K_{p+r+1}\),
implying that \( x \in \ZZ{r+1}{p} \).
Furthermore \( x \in K_{p+1} \), \( \partial_{r} x \in K_{p+r} \),
implying that \(\partial_{r-1} x \in K_{p+r}\) from which we conclude that
\( x \in \ZZ{r-1}{p+1} \).

Now if \( x \in \partial_{r} \ZZ{r}{p-r} \), then \( x = \partial_{r} y , y \in \ZZ{r}{p-r} \),
that is, \( x \in \ZZ{r}{p}\).
Therefore \( x \in K_p , \partial_{r+1} x=0 \), and it follows that
 \( x \in \ZZ{r+1}{p} \).

Since \( \ZZ{r+1}{p}\supset \partial_{r} \ZZ{r}{p-r} \), \( \ZZ{r+1}{p} \cap \ZZ{r-1}{p+1}
= \ZZ{r}{p+1} \),
one has
\begin{eqnarray}
H^p(\EE{r}{\cdot})&=&
\ZZ{r+1}{p}/(\partial_{r} \ZZ{r}{p-r}+\ZZ{r}{p+1})=\EE{r+1}{p} .
\end{eqnarray}
In this way we translate normal form problems into cohomology.
\end{proof}

\end{document}